\documentclass[12pt,leqno]{article}
\title{Some applications of Kummer and Stickelberger relations  }
\author{Roland Qu\^eme}
\usepackage{amsmath,amsthm,amstext,amsbsy}
\newtheorem{thm}{Theorem}[section]
\newtheorem{cor}[thm]{Corollary}

\newtheorem{lem}[thm]{Lemma}
\font\mathbb=msbm10
\newcommand{\N}{\mbox{\mathbb N}}

\newcommand{\Q}{\mbox{\mathbb Q}}
\newcommand{\Z}{\mbox{\mathbb Z}}

\newcommand{\modu}{\ \mbox{mod}\ }
\newcommand{\be}{\begin{equation}}
\newcommand{\ee}{\end{equation}}
\date{2006 april 19}
\begin{document}
\maketitle
\tableofcontents
\clearpage
\abstract
Roland Qu\^eme

13 avenue du ch\^ateau d'eau

31490 Brax

France

tel : 0561067020

cel : 0684728729

mailto: roland.queme@wanadoo.fr

home page: http://roland.queme.free.fr/

************************************

V10 - MSC Classification : 11R18;  11R29

************************************

Let $p$ be an odd  prime. Let ${\bf F}_p$ be the finite field of $p$ elements with no null part ${\bf F}_p^*$. Let $K_p=\Q(\zeta_p)$ be the $p$-cyclotomic field. 
Let $\pi$ be the prime ideal of $K_p$ lying over $p$.
Let $v$ be a primitive root $\modu p$.
In the sequel of this paper, for $n\in \Z$  let us note briefly  $v^n$ for $v^n \modu p$ with $1\leq v^n\leq p-1$. 
Let $\sigma :\zeta_p\rightarrow \zeta_p^v$ be a $\Q$-isomorphism of $K_p/\Q$.  
Let $G_p$ be the Galois group of $K_p/\Q$.
Let $P(\sigma)=\sum_{i=0}^{p-2}\sigma^i\times v^{-i},\quad P(\sigma)\in\Z[G_p]$.

We suppose that $p$ is an irregular prime. 
Let $C_p$ be the  $p$-class group of $K_p$. Let $\Gamma$ be a subgroup of $C_p$ of order $p$ annihilated by $\sigma-\mu$ with 
$\mu\in{\bf F}_p^*$.  
From Kummer, there exist    not principal prime ideals  $\mathbf q$  of $\Z[\zeta_p]$ of inertial degree $1$ with class $Cl(\mathbf q)\in \Gamma$.
Let $q$ be the prime number lying above $\mathbf q$.

Let $n$ be the smallest natural integer $1< n\leq p-2$ such that $\mu\equiv v^n\modu p$ for $\mu$ defined above. There exist     singular  numbers $A$ with $A \Z[\zeta_p]= \mathbf q^p$ and $\pi^{n}\ |\ A-a$  
where $a$ is a natural number. If $A$ is singular not primary then  $\pi^{n}\ \|\ A-a$ and if  
$A$ is singular  primary then  $\pi^{p}\ |\ A-a$. 
We prove,  by an application of Stickelberger relation to the prime ideal $\mathbf q$,   that now  we can {\it climb} up to the $\pi$-adic congruence:
\begin{enumerate}
\item
 $\pi^{2p-1} \ |\ A^{P(\sigma)}$ if $q\equiv 1\modu p$.
\item
$\pi^{2p-1} \ \|\ A^{P(\sigma)}$ if $q\equiv 1\modu p$ and $p^{(q-1)/p}\equiv 1\modu q$.
\item
 $\pi^{2p} \ |\ A^{P(\sigma)}$ if $q\not\equiv 1\modu p$.
\end{enumerate}
This property  of $\pi$-adic congruences on singular numbers is at the heart of this paper. 
\begin{enumerate}
\item
As a first example, in section \ref{s601193} p. \pageref{s601193} this $\pi$-adic improvement allows us to give an elementary straightforward proof that  the 
relative $p$-class group $C_p^-$ verifies  the following  congruence $\modu p$: 
with $v,m$ defined above, the congruence 
\begin{equation} 
\sum_{i=1}^{p-2} v^{(2m+1)(i-1)}\times(\frac{v^{-(i-1)}-v^{-i}\times v}{p})\equiv 0\modu p,
\end{equation}
is verified for $m$ taking $r^-$ different values $m_i,\quad i=1,\dots,r^-$ where $r^-$ is the rank of the  relative $p$-class group $C_p^-$
(result which can  also  be proved by annihilation of class group of $K_p$ by Stickelberger ideal $\in \Z[G_p]$).
A second example is a straightforward proof that if $\frac{p-1}{2}$ is odd then the Bernoulli Number $B_{(p+1)/2}\not\equiv 0\modu p$.
\item
The section \ref{s604191} p. \pageref{s604191} brings some results on connection between  singular primary numbers and  the stucture of the $p$-class group of $K_p$. 
\item
In the  last   section \ref{s604192} p. \pageref{s604192} we give  some  explicit  congruences  derived of Stickelberger  for  prime ideals $\mathbf q$ of inertial degree $f>1$.
\end{enumerate}

%

%
\section{Some definitions}\label{s601191}
In this section we give the definitions and  notations on cyclotomic fields,  $p$-class group,   singular  numbers,  primary and not primary,
used in this paper.
\begin{enumerate}
\item
Let $p$ be an odd prime. Let $\zeta_p$ be a root of the polynomial equation $X^{p-1}+X^{p-2}+\dots+X+1=0$.
Let $K_p$ be the $p$-cyclotomic field $K_p=\Q(\zeta_p)$. The ring of integers of $K_p$ is $\Z[\zeta_p]$.
Let $K_p^+$ be the maximal totally real subfield of $K_p$.
The ring of integers of $K_p^+$ is $\Z[\zeta_p+\zeta_p^{-1}]$ with group of units $\Z[\zeta_p+\zeta_p^{-1}]^*$.
Let $v$ be a primitive root $\modu p$ and $\sigma: \zeta_p\rightarrow \zeta_p^v$ be a $\Q$-isomorphism of $K_p$.
Let $G_p$ be the Galois group of $K_p/\Q$.
Let ${\bf F}_p$ be the finite field of cardinal $p$ with no null part  ${\bf F}_p^*$.
Let $\lambda=\zeta_p-1$. The prime ideal of $K_p$ lying over $p$ is $\pi=\lambda \Z[\zeta_p]$.
\item 
Suppose that $p$ is irregular. 
Let $C_p$ be the $p$-class group of $K_p$.
Let $r$ be the rank of $C_p$.
Let $C_p^+$ be the $p$-class group of $K_p^+$. Then $C_p=C_p^+\oplus C_p^-$ where $C_p^-$ is the relative $p$-class group.
\item
Let $\Gamma$ be a  subgroup of order $p$ of $C_p$ annihilated by $\sigma-\mu\in {\bf F}_p[G_p]$ with $\mu\in{\bf F}_p^*$. Then $\mu\equiv v^{n}\modu p$ with a natural integer $n,\quad 1< n\leq p-2$.
\item 
An  integer  $A\in \Z[\zeta_p]$ is said singular if $A^{1/p}\not\in K_p$ and if  there exists an   ideal $\mathbf a$  of $\Z[\zeta_p]$ such that 
$A \Z[\zeta_p]=\mathbf a^p$.
\begin{enumerate}
\item
\underline {If   $\Gamma\subset  C_p^-$: }
then there exists  singular integers $A$  with $A \Z[\zeta_p] =\mathbf a^p$ where $\mathbf a$ is a {\bf not} principal  ideal of $\Z[\zeta_p]$ 
verifying   simultaneously 
\begin{equation}\label{e512101}
\begin{split}
& Cl(\mathbf a)\in \Gamma,\\
& \sigma(A)=A^\mu\times\alpha^p,\quad \mu\in {\bf F}_p^*,\quad \alpha\in K_p,\\
&\mu\equiv v^{2m+1}\modu p, \quad m\in\N, \quad 1\leq m\leq \frac{p-3}{2},\\
&\pi^{2m+1} \ |\ A-a,\quad a\in\N,\quad 1\leq a\leq p-1,\\
\end{split}
\end{equation}
Moreover, this number $A$ verifies \begin{equation}\label{e512103}
A\times\overline{A}=D^p,
\end{equation}
for some integer $D\in O_{K_p^+}$.
\begin{enumerate}
\item 
This integer $A$ is singular not primary  if  $\pi^{2m+1} \ \|\ A-a$.
\item
This integer $A$ is singular  primary  if  $\pi^{p}\   |\ A-a^p$.
\end{enumerate}
\item
\underline {If   $\Gamma\subset C_p^+$: }
then there exists  singular integers $A$  with $A \Z[\zeta_p] =\mathbf a^p$ where $\mathbf a$ is a {\bf not} principal  ideal of $\Z[\zeta_p]$ 
verifying   simultaneously 
\begin{equation}\label{e6012210}
\begin{split}
& Cl(\mathbf a)\in \Gamma,\\
& \sigma(A)=A^\mu\times\alpha^p,\quad \mu\in {\bf F}_p^*,\quad \alpha\in K_p,\\
&\mu\equiv v^{2m}\modu p, \quad m\in\Z, \quad 1\leq m\leq \frac{p-3}{2},\\
&\pi^{2m} \ |\ A-a,\quad a\in\Z,\quad 1\leq a\leq p-1,\\
\end{split}
\end{equation}
Moreover, this number $A$ verifies \begin{equation}\label{e512103}
\frac{A}{\overline{A}}=D^p,
\end{equation}
for some number  $D\in K_p^+$.
\begin{enumerate}
\item 
This integer $A$ is singular not primary  if  $\pi^{2m} \ \|\ A-a$.
\item
This number $A$ is singular  primary  if  $\pi^{p}\  |\ A-a^p$.
\end{enumerate}
\end{enumerate}
\end{enumerate} 
%
\section{On Kummer and Stickelberger relation}\label{s601192}
\begin{enumerate}
\item
Here we fix a notation for the sequel. 
Let $v$ be a primitive root $\modu p$. For every integer $k\in\Z$   then  $v^k$ is understood $\modu p$ so $1\leq v^k\leq p-1$.
 If $k<0$ it is to be understood as $v^kv^{-k}\equiv 1\modu p$.

\item
Let $q\not=p$ be an odd prime. 
Let $\zeta_q$ be a root of the minimal polynomial equation $X^{q-1}+X^{q-2}+\dots+X+1=0$.
Let $K_q=\Q(\zeta_q)$ be the $q$-cyclotomic field.
The ring of integers of $K_q$ is $\Z[\zeta_q]$.
Here we fix a notation for the sequel. 
Let $u$ be a primitive root $\modu q$. For every integer $k\in\Z$   then  $u^k$ is understood $\modu q$ so $1\leq u^k\leq q-1$.
 If $k<0$ it is to be understood as $u^ku^{-k}\equiv 1\modu q$.
Let $K_{pq}=\Q(\zeta_p,\zeta_q)$. Then $K_{pq}$ is the compositum $K_pK_q$.
The ring of integers of $K_{pq}$ is $\Z[\zeta_{pq}]$.
\item
Let $\mathbf q$ be a prime ideal of $\Z[\zeta_p]$ lying over the prime $q$.
Let $m=N_{K_p/\Q}(\mathbf q)= q^f$ where $f$ is the smallest integer such that $q^f\equiv 1\modu p$.
If $\psi(\alpha)=a$ is the image of $\alpha\in \Z[\zeta_p]$ under the natural map 
$\psi: \Z[\zeta_p]\rightarrow \Z[\zeta_p]/\mathbf q$, then for 
$\psi(\alpha)=a\not\equiv 0$ define a character $\chi_{\mathbf q}^{(p)}$ on ${\bf F}_m=\Z[\zeta_p]/\mathbf q$ by
\begin{equation}
\chi_{\mathbf q}^{(p)}(a)=\{\frac{\alpha}{\mathbf q}\}_p^{-1}=\overline{\{\frac{\alpha}{\mathbf q}\}}_p,
\end{equation}
where $\{\frac{\alpha}{\mathbf q}\}=\zeta_p^c$ for some natural integer $c$,  is the $p^{th}$ power residue character $\modu \mathbf q$.
We define 
\begin{equation}\label{e6012211}
g(\mathbf q)=\sum_{x\in{\bf F}_m}(\chi_{\mathbf q}^{(p)}(x)\times 
\zeta_q^{Tr_{{\bf F}_m/{\bf F}_q}(x)})\in \Z[\zeta_{pq}],
\end{equation}
and $\mathbf G(\mathbf q)= g(\mathbf q)^p$.
It follows that $\mathbf G(\mathbf q)\in \Z[\zeta_{pq}]$.
Moreover $\mathbf G(\mathbf q)=g(\mathbf q)^p\in \Z[\zeta_p]$, see for instance Mollin \cite{mol} prop. 5.88 (c) p. 308
or Ireland-Rosen \cite{ire} prop. 14.3.1 (c) p. 208.
\end{enumerate}
%
The Stickelberger's relation is classically:
\begin{thm}\label{t12201}

In $\Z[\zeta_p]$ we have the ideal decomposition 
\begin{equation}\label{e512121}
\mathbf G(\mathbf q)\Z[\zeta_p]=\mathbf q^{S},
\end{equation}
with $S=\sum_{t=1}^{p-1} t\times \varpi_t^{-1}$
where  $\varpi_t\in Gal(K_p/\Q)$ is given by $\varpi_t: \zeta_p\rightarrow \zeta_p^t$.
\end{thm}
See for instance Mollin \cite{mol} thm. 5.109 p. 315 and Ireland-Rosen \cite{ire} thm. 2. p.209.
%
\subsection{On the structure of $\mathbf G(\mathbf q)$.}
In this subsection we are studying carefully the structure of $g(\mathbf q)$ and $\mathbf G(\mathbf q)$.
\begin{lem}\label{l512151}
If $q\not\equiv 1\modu p$ then $g(\mathbf q)\in \Z[\zeta_p]$.
\begin{proof}$ $
\begin{enumerate}
\item
Let $u$ be a primitive root $\modu q$. Let $\tau :\zeta_q\rightarrow \zeta_q^u$ be a $\Q$-isomorphism generating $Gal(K_q/\Q)$.
The isomorphism $\tau$ is extended to a $K_p$-isomorphism of $K_{pq}$ by 
$\tau:\zeta_q\rightarrow \zeta_q^u,\quad \zeta_p\rightarrow \zeta_p$.
Then  $g(\mathbf q)^p=\mathbf G(\mathbf q)\in \Z[\zeta_p]$ and so 
\begin{displaymath}
\tau(g(\mathbf q))^p=g(\mathbf q)^p,
\end{displaymath}
and it follows that there exists a natural integer $\rho$ with $\rho<p$ such that 
\begin{displaymath}
\tau(g(\mathbf q))= \zeta_p^\rho\times  g(\mathbf q).
\end{displaymath}
Then $N_{K_{pq}/K_p}(\tau(g(\mathbf q)))=\zeta_p^{(q-1)\rho}\times N_{K_{pq}/K_p}(g(\mathbf q))$ and so  $\zeta_p^{\rho(q-1)}=1$.
\item
If $q\not\equiv 1\modu p$, it implies that $\zeta_p^\rho=1$ and so that $\tau(g(\mathbf q))=g(\mathbf q)$ and thus that $g(\mathbf q)\in \Z[\zeta_p]$.
\end{enumerate}
\end{proof}
\end{lem}
%
Let us note in the sequel $g(\mathbf q)=\sum_{i=0}^{q-2} g_i\times \zeta_q^i$ with $g_i\in \Z[\zeta_p]$.
\begin{lem}\label{l512152}
If $q\equiv 1\modu p$ then $g_0=0$.
\begin{proof}
Suppose that $g_0\not=0$ and search for a contradiction:
we start of 
\begin{displaymath}
\tau(g(\mathbf q))= \zeta_p^\rho\times  g(\mathbf q).
\end{displaymath}
We have $g(\mathbf q)=\sum_{i=0}^{q-2} g_i\times \zeta_q^i$  and so 
$\tau(g(\mathbf q))=\sum_{i=0}^{q-2}  g_i\times \zeta_q^{i u}$, 
therefore 
\begin{displaymath}
\sum_{i=0}^{q-2} (\zeta_p^\rho g_i)\times \zeta_q^i=\sum_{i=0}^{q-2}  g_i\times \zeta_q^{i u},
\end{displaymath}
thus $g_0=\zeta_p^\rho \times g_0$ and so $\zeta_p^\rho=1$ which 
implies that $\tau(g(\mathbf q))=g(\mathbf q)$ and so $g(\mathbf q)\in \Z[\zeta_p]$.
Then $\mathbf G(\mathbf q)=g(\mathbf q)^p$ and so Stickelberger relation leads to 
$g(\mathbf q)^p \Z[\zeta_p] =\mathbf q^S$ where $S=\sum_{t=1}^{p-1}t\times\varpi_t^{-1}$.
Therefore $\varpi_1^{-1}(\mathbf q) \ \|\ \mathbf q^S$ because $q$ splits totally in $K_p/\Q$
 and $\varpi_t^{-1}(\mathbf q)\not=\varpi_{t^\prime}^{-1}(\mathbf q)$ for $t\not=t^\prime$. This case is not possible because the first member 
$g(\mathbf q)^p$ is a $p$-power.
\end{proof}
\end{lem}
%
Here we give an elementary computation of $g(\mathbf q)$ not involving directly the Gauss Sums.
\begin{lem}\label{l512152}
If $q\equiv 1\modu p$ then 
\begin{equation}\label{e512151}
\begin{split}
& \mathbf G(\mathbf q) = g(\mathbf q)^p,\\
&g(\mathbf q)=\zeta_q +\zeta_p^\rho\zeta_q^{u^{-1}}+\zeta_p^{2\rho}\zeta_q^{u^{-2}}+\dots \zeta_p^{(q-2)\rho}\zeta_q^{u^{-(q-2)}},\\
& g(\mathbf q)^p \Z[\zeta_p] =\mathbf q^S,\\
\end{split}
\end{equation}
for some natural number $\rho,\quad 1<\rho\leq p-1$.
\begin{proof} $ $
\begin{enumerate}
\item
We start of $\tau(g(\mathbf q))=\zeta_p^\rho\times g(\mathbf q)$ and so 
\begin{equation}\label{e512152}
\sum_{i=1} ^{q-2}g_i \zeta_q^{ui}=\zeta_p^\rho\times\sum_{i=1}^{q-2} g_i \zeta_q^i,
\end{equation}
which implies that $g_i=g_1\zeta_p^\rho $ for  $u\times i\equiv 1\modu q$ and so $g_{u^{-1}}=g_1\zeta_p^\rho$ (where $u^{-1}$ is to be understood by
$u^{-1}\modu q$,  so $1\leq u^{-1}\leq q-1)$.
\item
Then
$\tau^2(g(\mathbf q))=\tau(\zeta_p^{\rho} g(\mathbf q))=\zeta_p^{2\rho} g(\mathbf q)$.
Then 
\begin{displaymath}
\sum_{i=1} ^{q-2}g_i \zeta_q^{u^2i}=\zeta_p^{2\rho}\times (\sum_{i=1}^{q-2} g_i \zeta_q^i),
\end{displaymath}
which implies that $g_i=g_1\zeta_p^{2\rho}$ for $u^2\times i\equiv 1\modu q$ and so $g_{u^{-2}}=g_1\zeta_p^{2\rho}$.
\item
We continue up to 
$\tau^{(q-2)\rho}(g(\mathbf q))=\tau^{q-3}(\zeta_p^\rho g(\mathbf q))=\dots=\zeta_p^{(q-2)\rho} g(\mathbf q)$.
Then 
\begin{displaymath}
\sum_{i=1} ^{q-2}g_i \zeta_q^{u^{q-2}i}=\zeta_p^{(q-2)\rho}\times(\sum_{i=1}^{q-2} g_i \zeta_q^i),
\end{displaymath}
which implies that $g_i=g_1\zeta_p^{(q-2)\rho}$ for $u^{q-2}\times i\equiv 1\modu q$ and so $g_{u^{-(q-2)}}=g_1\zeta_p^{(q-2)\rho}$.
\item
Observe that $u$ is a primitive root $\modu q$ and so $u^{-1}$ is a primitive root $\modu q$.
Then it follows that 
$g(\mathbf q) =g_1\times (\zeta_q +\zeta_p^\rho\zeta_q^{u^{-1}}+\zeta_p^{2\rho}\zeta_q^{u^{-2}}+\dots \zeta_p^{(q-2)\rho}\zeta_q^{u^{-(q-2)}})$.
Let $U=\zeta_q +\zeta_p^\rho\zeta_q^{u^{-1}}+\zeta_p^{2\rho}\zeta_q^{u^{-2}}+\dots \zeta_p^{(q-2)\rho}\zeta_q^{u^{-(q-2)}}$.
\item
We prove now that $g_1\in \Z[\zeta_p]^*$. 
From Stickelberger relation $g_1^p \times U^p =\mathbf q^{S}$.
From $S=\sum_{i=1}^{p-1}\varpi_t^{-1}\times t$ it follows that 
$\varpi_t^{-1}(\mathbf q)^t\ \|\ \mathbf q^{S}$
and so that $g_1\not\equiv 0\modu \varpi_t^{-1}(\mathbf q)$
because $g_1^p$ is a $p$-power,
which implies that 
$g_1\in \Z[\zeta_p]^*$. 
Let us consider the relation(\ref{e6012211}). Let $x=1\in{\bf F}_q$, then $Tr_{{\bf F}_q/{\bf F}_q}(x)=1$ and $\chi_\mathbf q^{(p)}(1)=1^{(q-1)/p}\modu \mathbf q=1$
and thus the coefficient of $\zeta_q$ is $1$ and so $g_1=1$. 
\item
From Stickelberger,  $g(\mathbf q)^p \Z[\zeta_p]=\mathbf q^S$,
which achieves the proof.
\end{enumerate}
\end{proof}
\end{lem}
%
\paragraph{Remark:}
From
\begin{equation}
\begin{split}
&g(\mathbf q)=\zeta_q +\zeta_p^\rho\zeta_q^{u^{-1}}+\zeta_p^{2\rho}\zeta_q^{u^{-2}}+\dots +\zeta_p^{(q-2)\rho}\zeta_q^{u^{-(q-2)}},\\
&\Rightarrow \tau(g(\mathbf q))=\zeta_q^u +\zeta_p^\rho\zeta_q+\zeta_p^{2\rho}\zeta_q^{u^{-1}}+\dots +\zeta_p^{(q-2)\rho}\zeta_q^{u^{-(q-3)}},\\
&\Rightarrow\zeta^\rho\times g(\mathbf q)=\zeta^\rho\zeta_q +\zeta_p^{2\rho}\zeta_q^{u^{-1}}+\zeta_p^{3\rho}\zeta_q^{u^{-2}}+\dots +\zeta_p^{(q-1)\rho}
\zeta_q^{u^{-(q-2)}}\\
\end{split}
\end{equation}
and we can verify directly that $\tau(g(\mathbf q))=\zeta_p^\rho \times g(\mathbf q)$ for this expression of $g(\mathbf q)$, observing that $q-1\equiv 0\modu p$.
%
%
\begin{lem}\label{l12161}
Let $S=\sum_{t=1}^{p-1} \varpi_t^{-1}\times t$ where $\varpi_t$ is the $\Q$-isomorphism 
given by $\varpi_t:\zeta_p\rightarrow \zeta_p^t$ of $K_p$.
Let $v$ be a primitive root $\modu p$. Let $\sigma$ be the $\Q$-isomorphism of $K_p$ given by  $\zeta_p\rightarrow\zeta_p^v$.
Let $P(\sigma)=\sum_{i=0}^{p-2} \sigma^i\times v^{-i}\in\Z[G_p]$.
Then $S=P(\sigma)$.
\begin{proof}
Let us consider one term $\varpi_t^{-1} \times t$. Then $v^{-1}=v^{p-2}$ is a primitive root $\modu p$ because $p-2$ and $p-1$ are coprime and so there exists one and one $i$ such that 
$t=v^{-i}$. Then $\varpi_{v^{-i}}:\zeta_p\rightarrow \zeta_p^{v^{-i}}$ and so $\varpi_{v^{-i}}^{-1}:\zeta_p\rightarrow\zeta_p^{v^i}$
and so $\varpi_{v^{-i}}^{-1}=\sigma^i$ (observe that $\sigma^{p-1}\times v^{-(p-1)}=1$), which achieves the proof.
\end{proof}
\end{lem} 
%
\paragraph{Remark} : The previous lemma is  a verification of  the consistency of results in Ribenboim \cite{rib} p. 118, of Mollin \cite {mol} p. 315 and of Ireland-Rosen p. 209 with our computation. In the sequel we use Ribenboim notation more adequate for the factorization in ${\bf F}_p[G]$.
%
In that case the Stickelberger's relation is connected with the Kummer's relation on Jacobi resolvents, see for instance 
Ribenboim, \cite{rib} (2A) b. p. 118 and (2C) relation (2.6) p. 119.

\begin{lem}\label{l512162}$ $
If $q\equiv 1\modu p$ then 
\begin{enumerate}
\item
$g(\mathbf q)$ defined in relation (\ref{e512151}) is a  Jacobi resolvent: $g(\mathbf q)=<\zeta_p^\rho,\zeta_q>$.
\item
$\rho=-v$.
\end{enumerate}
\begin{proof}$ $
\begin{enumerate}
\item
Apply formula of Ribenboim \cite{rib} (2.2) p. 118 with 
$p=p, q=q, \zeta =\zeta_p,\quad \rho=\zeta_q, \quad n=\rho,\quad u=i,\quad m=1$ and $h=u^{-1}$
(where the left members notations $p, q, \zeta,\rho, n, u, m$ and $h$ are the Ribenboim notations).
\item
We start of 
$<\zeta_p^\rho,\zeta_q>=g(\mathbf q)$.
Then $v$ is a primitive root $\modu p$, so there exists  a natural integer $l$ such that 
$\rho \equiv v^l\modu p$.
By conjugation $\sigma^{-l}$ we get 
$<\zeta_p,\zeta_q>=g(\mathbf q)^{\sigma^{-l}}$.
Raising to $p$-power 
$<\zeta_p,\zeta_q>^p=g(\mathbf q)^{p\sigma^{-l}}$.
From lemma \ref{l12161} and Stickelberger relation
$<\zeta_p,\zeta_q>^p\Z[\zeta_p]=\mathbf q^{P(\sigma)\sigma^{-l}}$.
From Kummer's relation (2.6) p. 119 in Ribenboim \cite{rib}, we get 
$<\zeta_p,\zeta_q>^p\Z[\zeta_p]=\mathbf q^{P_1(\sigma)}$ with $P_1(\sigma)=\sum_{j=0}^{p-2}\sigma^j v^{(p-1)/2-j}$.
Therefore
$\sum_{i=0}^{p-2} \sigma^{i-l}v^{-i}=\sum_{j=0}^{p-2}\sigma^j v^{(p-1)/2-j}$.
Then 
$i-l\equiv j\modu p$ and $-i\equiv \frac{p-1}{2}-j\modu p$ (or  $i\equiv j-\frac{p-1}{2}\modu p$)
imply that 
$j-\frac{p-1}{2}-l\equiv j\modu p$,
so 
$l+\frac{p-1}{2}\equiv 0\modu p$,
so
$l\equiv -\frac{p-1}{2}\modu p$,
and 
$l\equiv \frac{p+1}{2}\modu p$,
thus 
$\rho\equiv v^{(p+1)/2}\modu p$ and finally $\rho= -v$.
\end{enumerate}
\end{proof}
\end{lem}
\paragraph{Remark}: The previous lemma allows to verify the consistency of our computation with Jacobi resultents used in Kummer (see Ribenboim 
p. 118-119).
%
\begin{lem}\label{l512161}
If $\mathbf q\equiv 1\modu p$ then $g(\mathbf q)\equiv -1\modu \pi$.
\begin{proof}
From $g(\mathbf q)=\zeta_q +\zeta_p^\rho\zeta_q^{u^{-1}}+\zeta_p^{2\rho}\zeta_q^{u^{-2}}+\dots +\zeta_p^{(q-2)\rho}\zeta_q^{u^{-(q-2)}}$, we see that 
$g(\mathbf q)\equiv \zeta_q +\zeta_q^{u^{-1}}+\zeta_q^{u^{-2}}+\dots +\zeta_q^{u^{-(q-2)}}\modu \pi$.
From $u^{-1}$ primitive root $\modu p$ it follows  that 
$1+\zeta_q +\zeta_q^{u^{-1}}+\zeta_q^{u^{-2}}+\dots +\zeta_q^{u^{-(q-2)}}=0$, which leads to the result.
\end{proof}
\end{lem}
%
%
It is possible to improve the previous   result  to:
\begin{lem}\label{l602061}$ $
Suppose that $q\equiv 1\modu p$.
If  $p^{(q-1)/p}\not\equiv 1\modu q$ then $\pi^p\ \|\ g(\mathbf q)^p+1$.
\begin{proof}

$ $
\begin{enumerate}
\item
We start of $g(\mathbf q)=\zeta_q +\zeta_p^\rho\zeta_q^{u^{-1}}+\zeta_p^{2\rho}\zeta_q^{u^{-2}}+\dots \zeta_p^{(q-2)\rho}\zeta_q^{u^{-(q-2)}}$,  
so
\begin{displaymath}
g(\mathbf q)=\zeta_q +((\zeta_p^\rho-1)+1)\zeta_q^{u^{-1}}+((\zeta_p^{2\rho}-1)+1)\zeta_q^{u^{-2}},
+\dots ((\zeta_p^{(q-2)\rho}-1)+1)\zeta_q^{u^{-(q-2)}}
\end{displaymath}
also
\begin{displaymath}
g(\mathbf q)=-1+(\zeta_p^\rho-1)\zeta_q^{u^{-1}}+(\zeta_p^{2\rho}-1)\zeta_q^{u^{-2}}
+\dots +(\zeta_p^{(q-2)\rho}-1)\zeta_q^{u^{-(q-2)}}.
\end{displaymath}
Then $\zeta_p^{i\rho}\equiv 1+i\rho\lambda\modu \pi^2$, so
\begin{displaymath} 
g(\mathbf q)\equiv -1+\lambda\times (\rho\zeta_q^{u^{-1}}+2\rho \zeta_q^{u^{-2}}
+\dots +(q-2)\rho)\zeta_q^{u^{-(q-2)}})\modu\lambda^2.
\end{displaymath}
Then 
$g(\mathbf q) =-1+\lambda U+\lambda^2V$
with 
$U=\rho\zeta_q^{u^{-1}}+2\rho \zeta_q^{u^{-2}}
+\dots +(q-2)\rho)\zeta_q^{u^{-(q-2)}}$ and $U,V\in\Z[\zeta_{pq}]$.
\item
Suppose that $\pi^{p+1} \ |\ g(\mathbf q)^p+1$ and search for a contradiction:
then, from $g(\mathbf q)^p =(-1+\lambda U+\lambda^2V)^p$,
it follows that 
 $p\lambda U+\lambda^pU^p\equiv 0\modu\pi^{p+1}$
and so $U^p-U\equiv 0\modu \pi$ because $p\lambda+\lambda^p\equiv 0\modu\pi^{p+1}$.
Therefore 
\begin{displaymath}
\begin{split}
& (\rho\zeta_q^{u^{-1}}+2\rho \zeta_q^{u^{-2}}
+\dots +(q-2)\rho)\zeta_q^{u^{-(q-2)}})^p-\\
&(\rho\zeta_q^{u^{-1}}+2\rho \zeta_q^{u^{-2}}
+\dots +(q-2)\rho)\zeta_q^{u^{-(q-2)}})\equiv 0\modu \lambda,\\
\end{split}
\end{displaymath}
and so 
\begin{displaymath}
\begin{split}
& (\rho\zeta_q^{pu^{-1}}+2\rho \zeta_q^{pu^{-2}}
+\dots +(q-2)\rho)\zeta_q^{pu^{-(q-2)}})\\
& -(\rho\zeta_q^{u^{-1}}+2\rho \zeta_q^{u^{-2}}
+\dots +(q-2)\rho)\zeta_q^{u^{-(q-2)}})\equiv 0\modu \lambda.\\
\end{split}
\end{displaymath}
\item
For any  natural  $j$ with $1\leq j\leq q-2$,  there must exist a natural $j^\prime$ with $1\leq j^\prime\leq q-2$ such that simultaneously: 
\begin{displaymath}
\begin{split}
& pu^{-j^\prime}\equiv u^{-j}\modu q\Rightarrow p\equiv u^{j^\prime-j}\modu q,\\
& \Rightarrow\rho  j^\prime \equiv \rho  j\modu \pi\Rightarrow  j^\prime-j\equiv 0 \modu p.\\
\end{split}
\end{displaymath}
Therefore $p\equiv u^{p\times \{(j^\prime-j)/p\}}\modu q$ and so $p^{(q-1)/p}\equiv u^{p\times (q-1)/p)\times \{(j^\prime-j)/p\}}\modu q$
thus $p^{(q-1)/p}\equiv 1\modu q$, contradiction.
\end{enumerate}
\end{proof}
\end{lem}
%
%
\subsection{A study of polynomial  $P(\sigma)=\sum_{i=0}^{p-2}\sigma^i v^{-i}$ of  $\Z[G_p]$.}
Recall that $P(\sigma)\in\Z[G_p]$ has been defined by $P(\sigma)=\sum_{i=0}^{p-2}\sigma^iv^{-i}$.
\begin{lem}\label{l512171}
\begin{equation}
P(\sigma)=\sum_{i=0}^{p-2}\sigma^i\times v^{-i}=v^{-(p-2)}\times \{\prod_{k=0,\ k\not=1}^{p-2}(\sigma-v^{k})\}+p\times R(\sigma),
\end{equation}
where $R(\sigma)\in\Z[G_p]$ with $deg(R(\sigma))<p-2$.
\begin{proof}
Let us consider the polynomial $R_0(\sigma)=P(\sigma)-v^{-(p-2)}\times \{\prod_{k=0,\ k\not=1}^{p-2}(\sigma-v^{k})\}$  in ${\bf F}_p[G_p]$.
Then $R_0(\sigma)$ is of degree smaller than $p-2$ and  the two 
polynomials $\sum_{i=0}^{p-2} \sigma^iv^{-i} $ and  $\prod_{k=0,\ k\not=1}^{p-2}(\sigma-v^{k})$ take a null value in ${\bf F}_p[G_p]$  when $\sigma$ takes  the $p-2$ different  values  $\sigma=v^k$ for $k=0,\dots, p-2,\quad k\not= 1$. Then $R_0(\sigma)=0$ in ${\bf F}_p[G_p]$ which leads to the result in 
$\Z[G_p]$.
\end{proof}
\end{lem}
Let us note in the sequel  
\begin{equation}\label{e12201}
T(\sigma)=v^{-(p-2)}\times\prod_{k=0,\ k\not=1}^{p-2}(\sigma-v^{k}).
\end{equation}
%
\begin{lem}\label{l512165}
\begin{equation}\label{e512172}
P(\sigma)\times (\sigma-v)= T(\sigma)\times(\sigma-v)+pR(\sigma)\times (\sigma-v)=p\times Q(\sigma),
\end{equation}
where $Q(\sigma)=\sum_{i=1}^{p-2}\delta_i\times \sigma^i\in\Z[G_p]$ is given by  
\begin{equation}
\begin{split}
& \delta_{p-2}= \frac{v^{-(p-3)}-v^{-(p-2)}v}{p},\\
& \delta_{p-3}= \frac{v^{-(p-4)}-v^{-(p-3)}v}{p},\\
& \vdots\\
& \delta_i=\frac{v^{-(i-1)}-v^{-i} v}{p},\\
&\vdots\\
& \delta_1 = \frac{1-v^{-1}v}{p},\\
\end{split}
\end{equation}
with $-p<\delta_i\leq 0$.
\begin{proof}
We start of  the relation in $\Z[G_p]$
\begin{displaymath}
P(\sigma)\times(\sigma-v)= v^{-(p-2)}\times \prod_{k=0}^{p-2} (\sigma-v^k)+p\times R(\sigma)\times(\sigma-v)=p\times Q(\sigma),
\end{displaymath}
with $Q(\sigma)\in\Z[G_p]$ because  $\prod_{k=0}^{p-2} (\sigma-v^k)=0$ in ${\bf F}_p[G_p]$ and so 
$\prod_{k=0}^{p-2} (\sigma-v^k)=p\times R_1(\sigma)$  in $\Z[G_p]$.
Then  we identify in $\Z[G_p]$ the  coefficients in the relation
\begin{displaymath}
\begin{split}
&(v^{-(p-2)}\sigma^{p-2}+v^{-(p-3)}\sigma^{p-3}+\dots+v^{-1}\sigma+1)\times(\sigma-v)=\\
&p\times (\delta_{p-2}\sigma^{p-2}+\delta_{p-3}\sigma^{p-3}+\dots+\delta_1\sigma+\delta_0),\\
\end{split}
\end{displaymath}
where $\sigma^{p-1}=1$.
\end{proof}
\end{lem}
\paragraph{Remark:}
\begin{enumerate}
\item
Observe that, with our notations,  $\delta_i\in \Z,\quad i=1,\dots,p-2$, but generally $\delta_i\not\equiv 0\modu p$. 
\item
We see also that $-p< \delta_i\leq 0$.
Observe also that $\delta_0=\frac{v^{-(p-2)}-v}{p}=0$.
\end{enumerate}
%
\begin{lem}\label{l601211}
The polynomial $Q(\sigma)$ verifies 
\begin{equation}\label{e6012110}
Q(\sigma)=\{(1-\sigma)(\sum_{i=0}^{(p-3)/2}\delta_i\times \sigma^i)+(1-v)\sigma^{(p-1)/2}\}\times(\sum_{i=0}^{(p-3)/2}\sigma^i).
\end{equation}
\begin{proof}
We start of $\delta_i=\frac{v^{-(i-1)}-v^{-i}v}{p}$.
Then  
\begin{displaymath}
\delta_{i+(p-1)/2}=\frac{v^{-(i+(p-1)/2-1)}-v^{-(i+(p-1)/2)}}{p}=
\frac{p-v^{-(i-1)}-(p-v^{-i})v}{p}=1-v-\delta_i.
\end{displaymath}
Then 
\begin{displaymath}
\begin{split}
& Q(\sigma)=\sum_{i=0}^{(p-3)/2} (\delta_i\times (\sigma^i-\sigma^{i+(p-1)/2}+(1-v)\sigma^{i+(p-1)/2})\\
& =(\sum_{i=0}^{(p-3)/2}\delta_i\times\sigma^i)\times (1-\sigma^{(p-1)/2)})+(1-v)\times \sigma^{(p-1)/2}\times (\sum_{i=0}^{(p-3)/2}\sigma^i),\\
\end{split}
\end{displaymath}
which leads to the result.
\end{proof}
\end{lem}

%
\subsection{$\pi$-adic congruences on the singular integers $A$}
From now we suppose  that the prime ideal $\mathbf q$ of $\Z[\zeta_p]$ has a class $Cl(\mathbf q)\in \Gamma$ where $\Gamma$ is  
a subgroup of order $p$ of $C_p$ previously defined,  with a singular integer $A$  given  by $A \Z[\zeta_p] = \mathbf q^p$.

In an other part, we know that the group of ideal classes of the cyclotomic field is generated by the ideal classes of prime ideals of degree $1$, see for instance 
Ribenboim, \cite{rib} (3A) p. 119.

%
\begin{lem}\label{l512163}$ $

$(\frac{g(\mathbf q}{\overline{g(\mathbf q})})^{p^2}=(\frac{A}{\overline{A}})^{P(\sigma)}$.
\begin{proof}
We start of 
$\mathbf G(\mathbf q)\Z[\zeta_p]= g(\mathbf q)^p \Z[\zeta_p]= \mathbf q^S$.
Raising to $p$-power we get 
$g(\mathbf q)^{p^2} \Z[\zeta_p]= \mathbf q^{pS}$.
But $A \Z[\zeta_p] = \mathbf q^p$, so 
\begin{equation}\label{e601071}
g(\mathbf q)^{p^2} \Z[\zeta_p]= A^{S}\Z[\zeta_p],
\end{equation}
so 
\begin{equation}\label{e601063}
 g(\mathbf q)^{p^2}\times  \zeta_p^w\times \eta= A^{S},\quad \eta\in \Z[\zeta_p+\zeta_p^{-1}]^*,
\end{equation} where $w$ is a natural number.
Therefore, by complex conjugation, we get 
$ \overline{g(\mathbf q)}^{p^2}\times\zeta_p^{-w}\times  \eta= \overline{A}^{S}.$
Then 
$ (\frac{g(\mathbf q)}{\overline{g(\mathbf q)}})^{p^2}\times\zeta_p^{2w}=(\frac{A}{\overline{A}})^{S}$.
From $A\equiv a\modu\pi^{2m+1}$ with $a$ natural integer, we get $\frac{A}{\overline{A}}\equiv 1\modu\pi^{2m+1}$  and so $w=0$.
Then $ (\frac{g(\mathbf q)}{\overline{g(\mathbf q)}})^{p^2}=(\frac{A}{\overline{A}})^{S}$.
\end{proof}
\end{lem}
\paragraph{Remark:} Observe that this lemma is true if either $q\equiv 1\modu p$ or $q\not\equiv 1\modu p$.
%
\begin{thm}\label{t601311}$ $
\begin{enumerate}
\item
$g(\mathbf q)^{p^2}=\pm A^{P(\sigma)}$.
\item
$g(\mathbf q)^{p(\sigma-1)(\sigma-v)}=\pm
(\frac{\overline{A}}{A})^{Q_1(\sigma)}$
where 
\begin{displaymath}
Q_1(\sigma)= (1-\sigma)\times (\sum_{i=0}^{(p-3)/2}\delta_i\times \sigma^i)+(1-v)\times \sigma^{(p-1)/2}.
\end{displaymath}
\end{enumerate}
\begin{proof}$ $
\begin{enumerate}
\item
We start of  $g(\mathbf q)^{p^2}\times\eta = A^{P(\sigma)}$ proved.
Then $g(\mathbf q)^{p^2(\sigma-1)(\sigma-v)}\times\eta^{(\sigma-1)(\sigma-v)}=
A^{P(\sigma)(\sigma-1)(\sigma-v)}$.
From lemma \ref{l601211}, we get 
\begin{displaymath}
P(\sigma)\times (\sigma-v)\times (\sigma-1)=p \times Q_1(\sigma)\times (\sigma^{(p-1)/2}-1),
\end{displaymath}
where 
\begin{displaymath}
Q_1(\sigma)= (1-\sigma)\times (\sum_{i=0}^{(p-3)/2}\delta_i\times \sigma^i)+(1-v)\times \sigma^{(p-1)/2}.
\end{displaymath}
Therefore
\begin{equation}\label{e602013}
g(\mathbf q)^{p^2(\sigma-1)(\sigma-v)}\times\eta^{(\sigma-1)(\sigma-v)}=
(\frac{\overline{A}}{A})^{p Q_1(\sigma)},
\end{equation}
and by conjugation
\begin{displaymath}
\overline{g(\mathbf q)}^{p^2(\sigma-1)(\sigma-v)}\times\eta^{(\sigma-1)(\sigma-v)}=
(\frac{A}{\overline{A}})^{p Q_1(\sigma)}.
\end{displaymath}
Multiplying these two relations we get, observing that $g(\mathbf q)\times \overline{g(\mathbf q)}=q^f$,
\begin{displaymath}
q^{f p^2(\sigma-1)(\sigma-v)}\times \eta^{2(\sigma-1)(\sigma-v)}=1,
\end{displaymath}
also 
\begin{displaymath}
\eta^{2(\sigma-1)(\sigma-v)}=1,
\end{displaymath}
and thus $\eta=\pm 1$ because $\eta\in\Z[\zeta_p+\zeta_p^{-1}]^*$, which with relation (\ref{e602013}) get $g(\mathbf q)^{p^2}=\pm A^{P(\sigma)}$ achieves the proof of the first part. 
\item
From relation (\ref{e602013}) we get 
\begin{equation}\label{e602014}
g(\mathbf q)^{p^2(\sigma-1)(\sigma-v)}=\pm 
(\frac{\overline{A}}{A})^{p Q_1(\sigma)},
\end{equation}
so 
\begin{equation}\label{e602031}
g(\mathbf q)^{p(\sigma-1)(\sigma-v)}=\pm \zeta_p^w\times 
(\frac{\overline{A}}{A})^{ Q_1(\sigma)},
\end{equation}
where $w$ is a natural number.
But $g(\mathbf q)^{\sigma-v}\in K_p$ and so $g(\mathbf q)^{p(\sigma-v)(\sigma-1)}\in (K_p)^p$, see for instance Ribenboim \cite{rib} (2A) b. p. 118.
and $(\frac{\overline{A}}{A})^{Q_1(\sigma)}\in (K_p)^p$ because $\sigma-\mu \ |\ Q_1(\sigma)$ in ${\bf F}_p[G_p]$ imply that $w=0$,
which achieves the proof of the second part.
\end{enumerate}
\end{proof}
\end{thm}
%
\paragraph{Remarks}
\begin{enumerate}
\item
Observe that this theorem is true either $q\equiv 1\modu p$ or $q\not\equiv 1\modu p$.
\item
$g(\mathbf q)\equiv -1\modu\pi$ implies that 
$g(\mathbf q)^{p^2}\equiv -1\modu \pi$.  
Observe that if $A\equiv a\modu\pi$ with $a$ natural number then 
$A^{P(\sigma)}\equiv a^{1+v^{-1}+\dots+v^{-(p-2)}}=a^{p(p-1)/2}\modu \pi\equiv \pm 1\modu \pi$ consistent with previous result.
\end{enumerate}
%
\begin{lem}\label{l60151}$ $
Let $q\not=p$ be an odd prime. Let $f$ be the smallest integer such that $q^f\equiv 1\modu p$. 
If $f$ is even then $g(\mathbf q)=\pm \zeta_p^wq^{f/2}$ for $w$ a natural number.
\begin{proof}$ $
\begin{enumerate}
\item
Let $\mathbf q$ be a prime ideal of $\Z[\zeta_p]$ lying over $q$. From $f$  even we get  
$\mathbf q=\overline{\mathbf q}$.
As in first section there exists singular numbers $A$ such that $A\Z[\zeta_p]=\mathbf q^p$. 
\item
From $\mathbf q=\overline{\mathbf q}$ we can choose  $A\in\Z[\zeta_p+\zeta_p^{-1}]$ and so $A=\overline{A}$.
\item
we have $g(\mathbf q)^{p^2}=\pm A^{P(\sigma)}$. 
From    lemma \ref{l512151} p. \pageref{l512151}, we know that $g(\mathbf q)\in \Z[\zeta_p]$.
\item
By complex conjugation $\overline{g(\mathbf q)^{p^2}}=\pm A^{P(\sigma)}$. 
Then  $g(\mathbf q)^{p^2}=\overline{g(\mathbf q)}^{p^2}$.
\item
Therefore $g(\mathbf q)^p =\zeta_p^{w_2}\times\overline{g(\mathbf q)}^p$ with 
$w_2$ natural number. As $g(\mathbf q)\in\Z[\zeta_p]$ this implies that $w_2=0$ and so $g(\mathbf q)^p =\overline{g(\mathbf q)}^p$.
Therefore  $g(\mathbf q) =\zeta_p^{w_3}\times\overline{g(\mathbf q)}$ with 
$w_3$ natural number.
But $g(\mathbf q)\times \overline{g(\mathbf q)}=q^f$ results of properties of power residue Gauss sums, see for instance 
Mollin prop 5.88 (b) p. 308.
Therefore $g(\mathbf q)^2=\zeta_p^{w_3}\times q^f$ and so $ (g(\mathbf q)\times\zeta_p^{-w_3/2})^2= q^f$ and thus 
$ g(\mathbf q)\times \zeta_p^{-w_3/2}=\pm  q^{f/2}$ wich achieves the proof.
\end{enumerate}
\end{proof}
\end{lem}  
%
\begin{thm}\label{l512164}$ $

\begin{enumerate}
\item
If $q\equiv 1 \modu p$ then $ A^{P(\sigma)}\equiv \delta \modu \pi^{2p-1}$ with $\delta\in\{-1,1\}$.
\item
If and only if $q\equiv 1 \modu p$ and $p^{(q-1)/p}\equiv 1\modu q$ then $\pi^{2p-1}\ \|\   A^{P(\sigma)} -\delta $ with $\delta\in\{-1,1\}$.
\item
If $q\not\equiv 1 \modu p$ then $ A^{P(\sigma)}\equiv \delta\modu \pi^{2p}$ with $\delta\in\{-1,1\}$.
\end{enumerate}
\begin{proof}$ $ 
\begin{enumerate}
\item
From lemma \ref{l512161}, we get $\pi^p\ |\ g(\mathbf q)^p+1$ and so $\pi^{2p-1}\ |\ g(\mathbf q)^{p^2}+1$. Then 
 apply theorem \ref{t601311}.
\item
Applying  lemma \ref{l602061} we get $\pi^p\ \|\ g(\mathbf q)^p+1$ and so $\pi^{2p-1}\ \|\ g(\mathbf q)^{p^2}+1$. Then 
apply theorem \ref{t601311}.
\item
From lemma \ref{l512151}, then $g(\mathbf q)\in\Z[\zeta_p]$ and so $\pi^{p+1}\ |\ g(\mathbf q)^p+1$ and also
$\pi^{2p}\ |\ g(\mathbf q)^{p^2}+1$.
\end{enumerate}
\end{proof}
\end{thm}
%
\paragraph{Remark:} 
If $C\in\Z[\zeta_p]$ is any semi-primary number  with  $C\equiv c\modu\pi^2$ with $c$ natural number 
we can only assert in general that $C^{P(\sigma)}\equiv  \pm 1\modu\pi^{p-1}$. For the singular numbers $A$ considered here we assert more: 
$A^{P(\sigma)}\equiv \pm 1\modu\pi^{2p-1}$.  We shall use this $\pi$-adic improvement in the sequel.

%

%
%
\section{Explicit polynomial congruences $\modu p$  connected to the $p$-class group}\label{s601193}
We deal of explicit polynomial congruences connected to  the $p$-class group when $p$ not divides the class number $h^+$ of $K_p^+$.
\begin{enumerate}
\item
We know  that the  relative    $p$-class group $C_p^-=\oplus_{k=1}^{r^-} \Gamma_k$ where $\Gamma_k$ are groups of order $p$ annihilated by $\sigma-\mu_k, 
\quad \mu_k \equiv v^{2m_k+1}\modu p,\quad 1\leq m_k\leq\frac{p-3}{2}$.
Let us consider the singular numbers $A_k,\quad k=1,\dots,r^-$, with $\pi^{2m_k+1} \ |\ A_k-\alpha_k$ with $\alpha_k$ natural number  defined  in  lemmas \ref{l512171} and \ref{l512165}.
From Kummer, the group of ideal classes of $K_p$ is generated by the classes of prime ideals of degree $1$ (see for instance Ribenboim \cite{rib} (3A) p. 119).
\item
In this section we shall explicit  a connection between  the  polynomial $Q(\sigma)\in\Z[G_p]$ and the structure of the relative $p$-class group $C_p^-$ of $K_p$.
\item
As another example we shall give an elementary proof in a straightforward way that if $\frac{p-1}{2}$ is odd then the Bernoulli Number $B_{(p+1)/2}\not\equiv 0\modu p$.
\end{enumerate}
%
\begin{thm}\label{t512171}
Let $p$ be an odd prime. Let $v$ be a primitive root $\modu p$. 
For $k=1,\dots,r^-$ rank of the $p$-class group of $K_p$  then 
\begin{equation}\label{e512191}
Q(v^{2m_k+1})=\sum_{i=1}^{p-2} v^{(2m_k+1)\times i}\times(\frac{v^{-(i-1)}-v^{-i}\times v}{p})\equiv 0\modu p,
\end{equation}
(or an other formulation  $\prod_{k=1}^{r^-}(\sigma-v^{2m_k+1})$ divides $Q(\sigma)$ in ${\bf F}_p[G_p]$).
\begin{proof}$ $
\begin{enumerate}
\item
Let us fix $A$ for one the singular numbers $A_k$ with $\pi^{2m+1}\  \|\ A-\alpha $ with $\alpha$ natural number  equivalent to 
$\pi^{2m+1}\ \|\  (\frac{A}{\overline{A}} -1)$,
equivalent to 
\begin{displaymath}
\frac{A}{\overline{A}} =1+\lambda^{2m+1}\times a,\quad a\in K_p, \quad v_{\pi}(a)=0.
\end{displaymath}
Then raising to $p$-power we get 
$(\frac{A}{\overline{A}})^p=(1+\lambda^{2m+1}\times a)^p\equiv 1+p\lambda^{2m+1} a\modu\pi^{p-1+2m+2}$ and so 
$\pi^{p-1+2m+1}\ \|\  (\frac{A}{\overline{A}})^p -1.$
\item
From theorem \ref{l512164} we get 
\begin{displaymath}
(\frac{A}{\overline{A}})^{P(\sigma)\times(\sigma-v)} =(\frac{A}{\overline{A}})^{pQ(\sigma)}\equiv 1\modu\pi^{2p-1}.
\end{displaymath}
We have shown that 
\begin{displaymath}
(\frac{A}{\overline{A}})^p=1+\lambda^{p-1+2m+1}b,\quad b\in K_p,\quad v_\pi(b)=0,
\end{displaymath}
then 
\begin{equation}\label{e512281}
(1+\lambda^{p-1+2m+1}b)^{Q(\sigma)}\equiv 1\modu  \pi^{2p-1}.
\end{equation}
\item
But $1+\lambda^{p-1+2m+1}b\equiv 1+p\lambda^{2m+1} b_1\modu \pi^{p-1+2m+2}$ with $b_1\in\Z,\quad b_1\not\equiv 0\modu p$. 
There exists a natural integer $n$ not divisible by $p$ such that 
\begin{displaymath}
(1+p\lambda^{2m+1} b_1)^n\equiv 1+p\lambda^{2m+1}\modu\pi^{p-1+2m+2}.
\end{displaymath}
Therefore 
\begin{equation}\label{e512192}
(1+p \lambda^{2m+1}b_1)^{n Q(\sigma)}\equiv (1+p \lambda^{2m+1})^{Q(\sigma)}\equiv 1\modu  \pi^{p-1+2m+2}.
\end{equation}
\item
Show that the possibility of climbing up the  step  $\modu \pi^{p-1+2m+2}$  implies that $\sigma-v^{2m+1}$ divides $Q(\sigma)$ in ${\bf F}_p[G_p]$:
we have $(1+p\lambda^{2m+1} )^\sigma=1+p\sigma(\lambda^{2m+1}) = 1+p(\zeta^v-1)^{2m+1} =1+p((\lambda+1)^v-1)^{2m+1}
\equiv 1+pv^{2m+1}\lambda^{2m+1}\modu\pi^{p-1+2m+2}.$
In an other part $(1+p\lambda^{2m+1})^{v^{2m+1}}\equiv 1+p v^{2m+1}\lambda^{2m+1}\modu\pi^{p-1+2m+2}$.
Therefore  
\begin{equation}\label{e512282}
(1+p\lambda^{2m+1})^{\sigma-v^{2m+1}}\equiv 1\modu\pi^{p-1+2m+2}.
\end{equation}
\item
By euclidean division of $Q(\sigma)$ by $\sigma-v^{2m+1}$ in ${\bf F}_p[G_p]$, we get 
\begin{displaymath}
Q(\sigma) = (\sigma-v^{2m+1}) Q_1(\sigma)+R
\end{displaymath}
with $R\in{\bf F}_p.$
From congruence (\ref{e512192}) and (\ref{e512282})  it follows that  $(1+p\lambda^{2m+1})^R\equiv 1\modu\pi^{p-1+2m+2}$ and so 
that $1+pR\lambda^{2m+1}\equiv 1\modu\pi^{p-1+2m+2}$ and finally 
that $R=0$. 
Then in ${\bf F}_p$ we have $Q(\sigma)=(\sigma-v^{2m+1})\times Q_1(\sigma)$ and so  $Q(v^{2m+1})\equiv 0\modu p$, or explicitly 
\begin{displaymath}
\begin{split}
& Q(v^{2m+1})=v^{(2m+1)(p-2)}\times \frac{v^{-(p-3)}-v^{-(p-2)}v}{p}+v^{(2m+1)(p-3)}\times\frac{v^{-(p-4)}-v^{-(p-3)}v}{p}+\dots\\
&+v^{2m+1}\times \frac{1-v^{-1}v}{p}\equiv 0\modu p,\\
\end{split}
\end{displaymath}
which achieves the proof.
\end{enumerate}
\end{proof}
\end{thm}

\paragraph{Remarks:}
\begin{enumerate}
\item
Observe   that it is the $\pi$-adic theorem  \ref{l512164} connected to Kummer-Stickelberger which allows to obtain this result.
\item
Observe that $\delta_i$ can also be written in the form $\delta_i=-[\frac{v^{-i}\times v}{p}]$ where $[x]$ is the integer part of $x$,
similar form also known in the literature.
\item
Observe that it is possible to get  other polynomials of $\Z[G_p]$ annihilating the relative $p$-class group $C_p^-$: for instance  
from  Kummer's formula on Jacobi cyclotomic functions we induce other  polynomials $Q_d(\sigma)$ annihilating the 
relative $p$-class group $C_p^-$  of $K_p$ : If $1\leq d\leq p-2$ define the set 
\begin{displaymath}
I_d=\{i\ |\ 0\leq i\leq p-2, \quad v^{(p-1)/2-i}+v^{(p-1)/2-i+ind_v(d)}> p\}
\end{displaymath}
where $ind_v(d)$ is the minimal integer $s$ such that $d\equiv v^s\modu p$.
Then the  polynomials $Q_d(\sigma)=\sum_{i\in I_d}\sigma^i$ for $d=1,\dots,p-2$ annihilate the $p$-class $C_p$ of $K_p$,
 see for instance Ribenboim \cite{rib}  relations (2.4) and (2.5) p. 119.
\item
See also in a more general context Washington, \cite{was} corollary 10.15 p. 198. 
\item
It is easy to verify the consistency of relation (\ref{e512191}) with the table of irregular primes and Bernoulli numbers in 
Washington, \cite{was} p. 410.
\end{enumerate}
%
An  immediate consequence is an explicit  criterium for $p$ to be a regular prime:
\begin{cor}\label{512301}
Let $p$ be an odd prime. Let $v$ be a primitive root $\modu p$.
If the congruence  
\begin{equation}\label{e512301}
\sum_{i=1}^{p-2} X^{i-1}\times(\frac{v^{-(i-1)}-v^{-i}\times v}{p})\equiv 0\modu p
\end{equation}
has no solution $X$ in $\Z$  with $X^{(p-1)/2}+1\equiv 0\modu p$ then the prime $p$ is regular.
\end{cor}
%
We give as another example a straightforward proof  of following lemma on Bernoulli Numbers  (compare elementary nature of this proof  with proof hinted by Washington in exercise 5.9 p. 85 using Siegel-Brauer theorem). 
\begin{lem}\label{l601011}
If $2m+1=\frac{p-1}{2}$ is odd then the Bernoulli Number $B_{(p+1)/2}\not\equiv 0\modu p$.
\begin{proof}
From previous corollary it follows that if $B_{(p+1)/2}\equiv 0\modu p$ implies that  $\sum_{i=1}^{p-2}v^{(2m+1)i}\times \delta^i\equiv 0\modu p$
where $2m+1=\frac{p-1}{2}$ because $v^{(p-1)/2}\equiv -1\modu p$. Then suppose that 
\begin{displaymath}
\sum_{i=1}^{p-2}(-1)^i\times (\frac{v^{-(i-1)}-v^{-i}\times v}{p})\equiv 0\modu p,
\end{displaymath}
and search for a contradiction:
multiplying by $p$
\begin{displaymath}
\sum_{i=1}^{p-2}(-1)^i\times (v^{-(i-1)}-v^{-i}\times v)\equiv 0\modu p^2,
\end{displaymath}
expanded to 
\begin{displaymath}
(-1+v^{-1}-v^{-2}+\dots-v^{-(p-3)})+(v^{-1}v-v^{-2}v+\dots+v^{-(p-2)}v)\equiv 0\modu p^2
\end{displaymath}
also 
\begin{displaymath}
(-1+v^{-1}-v^{-2}+\dots-v^{-(p-3)})+(v^{-1}-v^{-2}+\dots+v^{-(p-2)})v\equiv 0\modu p^2.
\end{displaymath}
Let us set $V=-1+v^{-1}-v^{-2}+\dots-v^{-(p-3)}+v^{-(p-2)}$.
Then we get 
$V-v^{-(p-2)}+v(V+1)\equiv 0\modu p^2$, and so 
$V(1+v)+v-v^{-(p-2)}\equiv 0\modu p^2$.
But $v=v^{-(p-2)}$ and so  $V\equiv 0\modu p^2$. 
But 
\begin{displaymath}
\begin{split}
& -V=1-v^{-1}+v^{-2}+\dots+v^{-(p-3)}-v^{-(p-2)}=S_1-S_2 \\
& S_1=1+v^{-2}+\dots+v^{-(p-3)},\\
& S_2=v^{-1}+v^{-3}+\dots+v^{-(p-2)}.\\
\end{split}
\end{displaymath}
$v^{-1}$ is a primitive root $\modu p$ and so $S_1+S_2=\frac{p(p-1)}{2}$.
Clearly $S_1\not= S_2$ because $\frac{p(p-1)}{2}$ is odd and so $-V=S_1-S_2\not=0$ and $-V\equiv 0\modu p^2$ with $|-V|<\frac{p(p-1)}{2}$, contradiction which achieves the proof.
\end{proof}
\end{lem}
%
\section{Singular primary numbers and Stickelberger relation}\label{s604191}
In this section we give some $\pi$-adic properties of singular numbers $A$ when they are primary.
Recall that $r, r^+, r^-$ are the ranks of the $p$-class groups $C_p, C_p^-,C_p^+$.
Recall that $C_p=\oplus_{i=1}^r \Gamma_i$ where $\Gamma_i$ are cyclic group of order $p$ annihilated by $\sigma-\mu_i$ with $\mu_i\in{\bf F}_p^*$.
\subsection {The case of $C_p^-$}
A classical result on structure of $p$-class group is that the relative $p$-class group $C_p^-$ is a direct sum $C_p^-=(\oplus_{i=1}^{r^+}\Gamma_i)\oplus(\oplus_{i=r^++1}^{r^-}\Gamma_i)$ where 
the subgroups $\Gamma_i,\quad i=1,\dots,r^+$ correspond to {\it singular primary} numbers  $A_i$ and 
where the subgroups $\Gamma_i,\quad i=r^++1,\dots,r^-$ corresponds to {\it singular not primary}  numbers $A_i$.
Let us fix one of these  singular primary numbers $A_i$ for $i=1,\dots,r^+$. Let $\mathbf q$ be a prime ideal of inertial degree $f$ such that 
$A\Z[\zeta_p]=\mathbf q^p$.
%
\begin{thm}\label{t601291}
Let $\mathbf q$ be a prime not principal ideal of $\Z[\zeta_p]$ of inertial degree $f$  with $Cl(\mathbf q)\in \Gamma\subset C_p^-$.
Suppose that  the prime number  $q$ above $\mathbf q$ verifies  $p\ \|\ q^f-1$ and 
that  $A$ is  a singular primary number with $A\Z[\zeta_p]=\mathbf q^p$.
Then 
\begin{equation}\label{e601291}
A\not\equiv 1 \modu\pi^{2p-1}.
\end{equation}
\begin{proof}$ $
\begin{enumerate}
\item
We start of the relation  
$g(\mathbf q)^{p^2} =\pm A^{P(\sigma)}$ proved  in theorem \ref{t601311}.
By conjugation we get 
$\overline{g(\mathbf q)}^{p^2} =\pm\overline{A}^{P(\sigma)}$.
Multiplying these two relations and observing that $g(\mathbf q)\times \overline{g(\mathbf q)}=q^f$ and $A\times\overline{A}=D^p$ with
$D\in\Z[\zeta_p+\zeta_p^{-1}]$ we get 
$q^{f p^2}= D^{pP(\sigma)}$, so $q^{f p}= D^{P(\sigma)}$ because $q, D\in\Z[\zeta_p+\zeta_p^{-1}]$ 
and, multiplying the exponent by $\sigma-v$,  we get 
$q^{f p(\sigma-v)} = D^{P(\sigma)(\sigma-v)}$ 
so $q^{f p(1-v)}= D^{p Q(\sigma)}$ from lemma \ref{l512165} p. \pageref{l512165} 
and thus 
\begin{equation}\label{e601301}
q^{f(1-v)}= D^{Q(\sigma)}.
\end{equation}
\item
Suppose that 
$\pi^{2p-1}\ |\ A-1$.
Then $\pi^{2p-1}\ |\ \overline{A}-1$, so 
$\pi^{2p-1}\ |\ D^{p}-1$
and so 
$\pi^{p}\ |\ D-1$ and so $\pi^{p}\ |\ D^{ Q(\sigma)}-1$, 
thus $\pi^{p}\ |\ q^{f(1-v)}-1$ and finally $\pi^{p}\ |\ q^f-1$, contradiction with $\pi^{p-1}\ \|\ q^f-1$.
\end{enumerate}
\end{proof}
\end{thm} 
%
In the following theorem we obtain a result of same nature  which can be applied generally to a wider range of singular primary numbers $A$
if we assume  simultaneously the two hypotheses $q\equiv 1\modu p$ and $p^{(q-1)/p}\equiv 1\modu q$. 
%
\begin{thm}\label{t601291}
Let $\mathbf q$ be a prime not principal ideal of $\Z[\zeta_p]$ of inertial degree $f=1$  with $Cl(\mathbf q)\in \Gamma\subset C_p$.
Let $A$ be  a singular primary number with $A\Z[\zeta_p]=\mathbf q^p$.
If $p^{(q-1)/p}\equiv 1\modu q$ then there exists no natural integer $a$ such that  
\begin{equation}\label{e601291}
A\equiv a^p \modu\pi^{2p}.
\end{equation}
\begin{proof}$ $
Suppose that $A\equiv a^p\modu\pi^{2p}$ and search for a contradiction.
We start of relation  $g(\mathbf q)^{p^2}=\pm A^{P(\sigma)}$ proved in theorem \ref{t601311} p. \pageref{t601311}.
Therefore 
$g(\mathbf q)^{p^2}\equiv \pm a^{p P(\sigma)}\modu\pi^{2p}$,
so
\begin{displaymath}
g(\mathbf q)^{p^2}\equiv \pm a^{p( v^{-(p-2)}+\dots+v^{-1}+1)}\modu\pi^{2p},
\end{displaymath}
so 
\begin{displaymath}
g(\mathbf q)^{p^2}\equiv \pm a^{p^2(p-1)/2}\modu\pi^{2p}.
\end{displaymath}
But $a^{p^2(p-1)/2}\equiv \pm 1\modu\pi^{2p}$.
It should imply that $g(\mathbf q)^{p^2}\equiv\pm 1 \modu \pi^{2p}$, so that $g(\mathbf q)^p\equiv \pm 1\modu \pi^{p+1}$ which 
contradicts lemma \ref{l602061} p. \pageref{l602061}.
\end{proof}
\end{thm}
%
\subsection{On principal prime ideals of $K_p$ and Stickelberger relation}
The Stickelberger relation and its consequences on prime ideals $\mathbf q$ of $\Z[\zeta_p]$  is meaningful even if $\mathbf q$ is a principal ideal.
\begin{thm}\label{t602061}
Let $q_1\in\Z[\zeta_p]$ with $q_1\equiv a\modu \pi^{p+1}$ where $a\in\Z,\quad  a\not\equiv 0\modu p$.
If $q=N_{K_p/\Q}(q_1)$ is a prime number then $p^{(q-1)/p}\equiv 1\modu q$.
\begin{proof}
From Stickelberger relation
$g(q_1\Z[\zeta_p])^p\Z[\zeta_p]=q_1^{P(\sigma)}\Z[\zeta_p]$ and so there exists $\varepsilon\in\Z[\zeta_p]^*$ such that 
$g(q_1\Z[\zeta_p])^p=q_1^{P(\sigma)}\times \varepsilon$ and so  
\begin{displaymath}
(\frac{g(q_1\Z[\zeta_p])}{\overline{g(q_1\Z[\zeta_p]}})^p= (\frac{q_1}{\overline{q_1}})^{P(\sigma)}.
\end{displaymath}
From hypothesis $\frac{q_1}{\overline{q_1}}\equiv 1\modu\pi^{p+1}$ and so 
$(\frac{g(q_1\Z[\zeta_p]}{\overline{g(q_1\Z[\zeta_p]}})^p\equiv 1\modu \pi^{p+1}$. From lemma \ref{l602061} p. \pageref{l602061} it follows that
$p^{(q-1)/p}\equiv 1\modu q$.
\end{proof}
\end{thm}
%
%
\section{Stickelberger's relation for prime ideals $\mathbf q$ of  inertial degree $f>1$.}\label{s604192}
Recall that the Stickelberger's relation is $g(\mathbf q)^p=\mathbf q^{S}$
where $S=\sum_{i=0}^{p-2}\sigma^iv^{-i}\in\Z[G_p]$.
We apply  Stickelberger's relation with the same method  to prime ideals $\mathbf q$ of  inertial degree $f>1$. 
Observe, from lemma \ref{l512151} p.
\pageref{l512151},  that  $f>1$ implies    $g(\mathbf q)\in\Z[\zeta_p]$. 

\paragraph{A definition:} we say that the prime ideal $\mathbf c$ of  a number field $M$ is $p$-principal if the component of the class group
$<Cl(\mathbf c)>$ in $p$-class group $D_p$ of $M$ is trivial. 
%
%
%

\begin{lem}\label{l603302}
Let $p$ be an odd prime. Let $v$ be a primitive root $\modu p$.
Let $q$ be an odd prime with $q\not=p$. Let $f$ be the smallest integer such that $q^f\equiv 1\modu p$ and 
let $m=\frac{p-1}{f}$. 
Let $\mathbf q$ be an  prime   ideal of $\Z[\zeta_p]$ lying over $q$. 
If $f>1$ then $g(\mathbf q)\in \Z[\zeta_p]$  and $g(\mathbf q)\Z[\zeta_p]= \mathbf q^{S_2}$ 
where 
\begin{equation}\label{e604013}
S_2=\sum_{i=0}^{m-1}(\frac{\sum_{j=0}^{f-1} v^{-(i+jm)}}{p})\times\sigma^i\in\Z[G_p].
\end{equation}
\begin{proof}$ $
\begin{enumerate}
\item
Let $p=fm+1$. Then 
$N_{K_p/\Q}(\mathbf q)=q^f$
and 
$\mathbf q=\mathbf q^{\sigma^m}=\dots=\mathbf q^{\sigma^{(f-1)m}}$.
The sum $S$ defined in lemma \ref{l12161} p.\pageref{l12161} can be written
\begin{displaymath}
S=\sum_{i=0}^{m-1}\sum_{j=0}^{f-1}\sigma^{i+j m}v^{-(i+jm)}.
\end{displaymath}
\item
From Stickelberger's relation seen in theorem \ref{t12201} p. \pageref{t12201}, 
then $g(\mathbf q)^p \Z[\zeta_p]=\mathbf q^{S}$.
Observe that, from hypothesis, $\mathbf q=\mathbf q^{\sigma^m}=\dots=\mathbf q^{\sigma^{(f-1)m}}$ so Stickelberger's relation implies that 
$g(\mathbf q)^p \Z[\zeta_p]=\mathbf q^{S_1}$
with
\begin{displaymath}
S_1=\sum_{i=0}^{m-1}\sum_{j=0}^{f-1}\sigma^{i}v^{-(i+jm)} = p\times \sum_{i=0}^{m-1} (\frac{\sum_{j=0}^{f-1} v^{-(i+jm)}}{p})\times  \sigma^i,
\end{displaymath}
where $(\sum_{i=0}^{f-1} v^{-(i+jm)})/p\in\Z$ because $v^{-m}-1\not\equiv 0\modu p$.
\item
Let $S_2=\frac{S_1}{p}$.
From lemma \ref{l512151} p.
\pageref{l512151} we know that  $f>1$ implies that   $g(\mathbf q)\in\Z[\zeta_p]$. 
Therefore 
\begin{displaymath}
g(\mathbf q) \Z[\zeta_p]=\mathbf q^{S_2}, \quad g(\mathbf q)\in\Z[\zeta_p].
\end{displaymath}
\end{enumerate}
\end{proof}
\end{lem}

%
%
It is possible to derive some explicit congruences in $\Z$ from this result.
\begin{lem}\label{l604012}
Let $p$ be an odd prime. Let $v$ be a primitive root $\modu p$.
Let $q$ be an odd prime with $q\not=p$. Let $f$ be the smallest integer such that $q^f\equiv 1\modu p$ and 
let $m=\frac{p-1}{f}$. 
Let $\mathbf q$ be an  prime   ideal of $\Z[\zeta_p]$ lying over $q$. 
\begin{enumerate}
\item
If $f>1$ and if $\mathbf q$ is not $p$-principal ideal  there exists  a natural integer $l,\ 1\leq l< m$ such that 
\begin{equation}\label{e604015}
\sum_{i=0}^{m-1}(\frac{\sum_{j=0}^{f-1} v^{-(i+jm)}}{p})\times v^{lfi}\equiv 0\modu p,
\end{equation}
\item
If for all natural integers $l$  such that $1\leq l<m$ 
\begin{equation}\label{e604011}
\sum_{i=0}^{m-1}(\frac{\sum_{j=0}^{f-1} v^{-(i+jm)}}{p})\times v^{lfi}\not\equiv 0\modu p,
\end{equation}
then $\mathbf q$ is $p$-principal
\end{enumerate}
\begin{proof}$ $
\begin{enumerate}
\item
Suppose that $\mathbf q$ is not $p$-principal.
Observe at first that congruence (\ref{e604015}) with $l=m$ should imply that  
$\sum_{i=0}^{m-1}(\sum_{j=0}^{f-1} v^{-(i+jm)})/p)\equiv 0\modu p$ or 
$\sum_{i=0}^{m-1}\sum_{j=0}^{f-1} v^{-(i+jm)}\equiv 0\modu p^2$  which is not possible because 
$v^{-(i+jm)}=v^{-(i^\prime+ j^\prime m)}$ implies that $j=j^\prime$ and $i=i^\prime$ and so that 
$\sum_{i=0}^{m-1}\sum_{j=0}^{f-1} v^{-(i+jm)}=\frac{p(p-1)}{2}$.
\item
The polynomial $S_2$ of lemma \ref{l603302} annihilates the 
not $p$-principal ideal $\mathbf q$  
in ${\bf F}_p[G_p]$  only if there exists $\sigma-v^n$ dividing $S_2$ in ${\bf F}_p[G_p]$. From $\mathbf q^{\sigma^m-1}=1$ it follows also that 
$\sigma-v^n\ |\ \sigma^m-1$. But $\sigma-v^n\ |\ \sigma^m-v^{nm}$ and so $\sigma-v^n\ |\ v^{nm}-1$, thus $nm\equiv 0\modu p-1$,  so
$n\equiv 0 \modu f$ and $n=lf$.
Therefore if $\mathbf q$ is not $p$-principal there exists  a natural integer $l,\ 1\leq l< m$ such that 
\begin{equation}\label{e604014}
\sum_{i=0}^{m-1}(\frac{\sum_{j=0}^{f-1} v^{-(i+jm)}}{p})\times v^{lfi}\equiv 0\modu p,
\end{equation}
\item
The relation (\ref{e604011}) is an imediate consequence of previous part of the proof.
\end{enumerate}
\end{proof}
\end{lem}

%
As an example we deal with  the case $f=\frac{p-1}{2}$.
\begin{cor}\label{c604011}
If $p\equiv 3\modu 4$ and if $f=\frac{p-1}{2}$ then $\mathbf q$ is $p$-principal. 
\begin{proof}
We have $f=\frac{p-1}{2}$, $m=2$ and $l=1$.
Then 
\begin{equation}\label{e604012}
\Sigma=\sum_{i=0}^{m-1}(\frac{\sum_{j=0}^{f-1} v^{-(i+jm)}}{p})\times v^{lfi}=
\frac{\sum_{j=0}^{(p-3)/2} v^{-2j}}{p}-\frac{\sum_{j=0}^{(p-3)/2} v^{-(1+2j)}}{p}.
\end{equation}
$\Sigma\equiv 0\modu p$ should imply that
$\sum_{j=0}^{(p-3)/2} v^{-2j}-\sum_{j=0}^{(p-3)/2} v^{-(1+2j)}\equiv 0\modu p^2$.
But $\sum_{j=0}^{(p-3)/2} v^{-2j}+\sum_{j=0}^{(p-3)/2} v^{-(1+2j)}=\frac{p(p-1)}{2}$ is odd, which achieves the proof.
\end{proof}
\end{cor}

%

%

\clearpage
\begin{thebibliography}{9}
\bibitem{ire} K. Ireland, M. Rosen, \textit{A Classical Introduction to Modern Number Theory}, Springer-Verlag, 1982.
\bibitem{koc} H. Koch, \textit{Algebraic Number Theory}, Springer, 1997.
\bibitem{mol} R.A. Mollin, \textit{Algebraic Number Theory}, Chapman  and Hall/CRC, 1999.
\bibitem{nar} W. Narkiewicz, \textit{Elementary and Analytic Theory of Numbers}, Springer-verlag, 1990.
\bibitem{rib} P. Ribenboim, \textit{13 Lectures on Fermat's Last Theorem}, Springer-Verlag, 1979.
\bibitem{ri2} P. Ribenboim, \textit{Classical Theory of Algebraic Numbers}, Springer, 2001.
\bibitem{was} L.C. Washington, \textit{Introduction to cyclotomic fields, second edition}, Springer, 1997.
\end{thebibliography}
\end{document}